\documentclass[11pt]{article}
\usepackage{latexsym}
\usepackage{theorem}
\usepackage{graphicx}
\usepackage{amsmath,color}
\usepackage{amsfonts}
\usepackage{natbib}
\usepackage{soul}

\theorembodyfont{\rmfamily}
\newtheorem{theorem}{Theorem}

\newtheorem{lemma}[theorem]{Lemma}

\newtheorem{definition}[theorem]{Definition}

\theoremstyle{break}

\begin{document}

\title{A Game Theoretic Approach to a Problem in \\ Polymatroid Maximization}

\date{}
\author{Lisa Hellerstein\thanks{Department of Computer Science and Engineering, New York University Tandon School of Engineering, New York, New York 11201, lisa.hellerstein@nyu.edu} \and Thomas Lidbetter\thanks{Department of Management Science and Information Systems, Rutgers Business School, Newark, NJ 07102, tlidbetter@business.rutgers.edu %(corresponding author)
}}

\providecommand{\keywords}[1]{\textbf{\textbf{Keywords:}} #1}

%\newpage

\maketitle

\begin{abstract}
\noindent We consider the problem of maximizing the minimum (weighted) value of all components of a vector over a polymatroid. This is a special case of the lexicographically optimal base problem introduced and solved by Fujishige. We give an alternative formulation of the problem as a zero-sum game between a maximizing player whose mixed strategy set is the base of the polymatroid and a minimizing player whose mixed strategy set is a simplex. We show that this game and three variations of it unify several problems in search, sequential testing and queuing. We give a new, short derivation of optimal strategies for both players and an expression for the value of the game. Furthermore, we give a characterization of the set of optimal strategies for the minimizing player and we consider special cases for which optimal strategies can be found particularly easily.
\end{abstract}

\keywords{Game theory; search games; sequential testing; queuing}
 %\newpage
\section{Introduction}
\label{sec:introduction}

A well understood problem in combinatorial optimization is that of maximizing a linear function over a polymatroid. As shown in \cite{edmonds1970submodular}, the solution of the problem is given by a simple greedy algorithm whose output is some vertex of the  base of the polymatroid. A similar algorithm can be used to minimize a linear function over a contrapolymatroid.  (All concepts will be defined precisely in Section~\ref{sec:main}).

Many optimization problems can be viewed as a special case of this problem. The general approach is to associate some ``performance vector'' with each possible choice of feasible solution to the problem in question, then to show that the convex hull $\mathcal B$ of these vectors is the base of a polymatroid or a contrapolymatroid. The objective function is then expressed as a linear function over $\mathcal B$, so that it can be optimized using the classic greedy algorithm.  

One example of such a problem is the single machine scheduling problem $1||\sum w_j C_j$ of choosing what order to process a finite set of jobs with given processing times to minimize their weighted sum of completion times: see \cite{queyranne1993structure} and also \cite{queyranne1994polyhedral}. \cite{agnetis2009sequencing} showed that another scheduling problem, introduced by \cite{stadje1995selecting}, in which an unreliable machine sequentially processes a set of jobs, can similarly be solved by maximizing a linear function over a polymatroid. \cite{kodialam2001throughput} had previously studied this same polymatroid to solve a different, but related problem in sequential testing. By considering so-called conservation laws, \cite{federgruen1988characterization} showed that the performance space of several multiclass queueing systems have a polymatroid structure, and this was extended to many other queueing problems by \cite{shanthikumar1992multiclass}. %Some cases of the well-known $c \mu$ priority rule for queueing can be derived by considering a linear optimization problem over the base of a polymatroid.

%\cite{stadje1995selecting} introduced a single machine scheduling problem where a finite set of jobs each has a given reward and a given probability of being successfully processed. If the machine fails to process a job it cannot process any further jobs.  The problem is to choose in which order to schedule the jobs to maximize the total expected reward.  Like the first scheduling problem described above, this problem is also solved by a simple index rule. More recently, the same problem was studied independently under the nomenclature of the {\em unreliable jobs problem} by \cite{agnetis2009sequencing}, who considered performance vectors corresponding to the probabilities that the jobs are successfully processed under a given schedule, and showed the convex hull of these vectors is the base of a polymatroid.  Hence, by writing the objective of the problem as a linear function over this polyhedron, the index rule for this problem can also be derived from the greedy algorithm for optimizing a linear function over a polymatroid. 

%More recently, \cite{agnetis2009sequencing} introduced a new single machine scheduling problem called the unreliable jobs problem, and showed that its solution (for the case of one machine) can also be derived from the greedy algorithm for maximizing a linear function over a polymatroid.

In this paper we focus on a max-min version of the classic problem of maximizing a linear function over a polymatroid. This max-min problem is a special case of the {\em lexicographically optimal base} problem, introduced and solved by~\cite{fujishige1980lexicographically}. We equivalently view the max-min problem as a zero-sum game between a maximizer whose pure strategies are the vertices of the base of an $n$-dimensional polymatroid and a minimizer whose pure strategies are the $n$ coordinate directions. Although the problem is already solved, we give a new, concise derivation of the solution using our game theoretic approach. 

More importantly, we show that our problem provides a unifying framework for many problems in search games, sequential testing and queueing; some known and some new.  Search games are two-person zero-sum games, where one player hides a ``target'' which the other player must locate. (See \cite{alpern2003theory} or \cite{hohzaki2016search} for an overview on the search games literature.)  In this paper we solve a case of the {\em weighted search game} introduced by \cite{yolmeh-gursoy}, where a Searcher aims to minimize a weighted time to find a target hidden among a finite number of locations with varying weights and search times.  We extend the weighted search game to incorporate the variable speed search paradigm of \cite{alpern2014searching}, and give a solution to this problem too. We show that the solution of a search and rescue game introduced by \cite{lidbetter2020search} also follows from a corollary of our main results; furthermore we solve a more elaborate search and rescue game. 

%LH I edited this slightly.  
We show that our approach yields an alternative solution to a problem in sequential testing previously solved by \cite{kodialam2001throughput} and \cite{condon2009algorithms}, in which operators sequentially perform tests on some tuples until obtaining a negative test, and the objective is to find a randomized routing of tuples to maximize throughput.

We also point out that our main problem can be used to address some max-min (or min-max) multiclass queueing problems, which, as far as we know, have not previously been considered in the literature. Although there are several possible applications, we consider one concrete example of a multiclass queueing problem in which one server processes jobs with exponentially distributed service times that arrive according to a Poisson process. The objective is to choose a randomized priority rule to minimize the maximum expected holding cost of any job class in the steady state of the system. This problem is a special case of our main problem.

% consider a max-min problem over the base $\mathcal{B}$ of a polymatroid. Rather than seeking to maximize a linear function $w^T x$ over $x \in \mathcal B$ we seek to maximize $\min_i w_i x_i$ over $x \in \mathcal B$. An equivalent approach that we prefer to take is to consider the problem as a zero-sum game between a maximizer whose pure strategies are the set of vertices $x$ of $\mathcal{B}$ and a minimizer whose pure strategies are the coordinates $i$, where the payoff of the game is $w_i x_i$. To the best of the author's knowledge, this problem has not been considered before in the literature.

In addition, we consider special cases of our main problem, where the payoff function satisfies certain monotonicity properties that we define later.  Although these special cases are more limited than the main problem, they include a number of particular problems, previously studied in the search games and sequential testing literature, which admit simpler solutions than the solutions to the main problem.

In Section~\ref{sec:main} we review the notion of a polymatroid and the classic greedy algorithm of \cite{edmonds1970submodular}. We then describe our main problem, framing it as a zero-sum game between a maximizer whose pure strategies are the set of vertices $x$ of the base of a polymatroid and a minimizer whose pure strategies are the coordinates $i$. The payoff of the game is $w_i x_i$ for some fixed positive weights $w$ (in contrast to the classic problem of \cite{edmonds1970submodular} where the objective is $w^T x$). We also describe three variations of the game involving contrapolymatroids and min-max objectives.

In Section~\ref{sec:apps}, we show that our problem and its variations unify several search games. We also make a link to a sequential testing problem and discuss further special cases of our problems in the field of queueing theory.

We give optimal strategies for both players in the main version of our game and an expression for its value in Section~\ref{sec:solution}. Our game theoretic angle on the problem yields insights that were not captured in the work of \cite{fujishige1980lexicographically} on the more general problem. We give a complete characterization of the set of optimal strategies for Player~2.  We also define special cases of the problem where the payoff function of the game satisfies certain monotonicity properties. 
While the value of the game and optimal strategies for both players can always be found in strongly polynomial time (in the dimension of the polymatroid), we show that in these special cases, the value of the game can be found particularly quickly.

We discuss the three variations of our game in Section~\ref{sec:min}, and in Section~\ref{sec:impl} we detail the implications of our results on the applications described in Section~\ref{sec:apps}. 
Finally, in Section~\ref{sec:monotone}, we further consider one of the special cases of our game where the payoff function satisfies a monotonicity condition.  
For this case, we give an efficient procedure that implements an optimal strategy for Player~1. 
The support of this strategy is of exponential size, and the procedure does not output an explicit representation of it as a convex combination of pure strategies.  Instead, the procedure can be used to efficiently generate a
pure strategy, drawn from the support of this optimal strategy with the appropriate probability. 

The results we give in Sections~\ref{sec:solution} and~\ref{sec:monotone} for the special cases are inspired by and generalize results in \cite{condon2009algorithms} for particular sequential testing problems.

\section{Problem Statement} \label{sec:main}

In this section we define and solve our main problem, then describe some applications.

\subsection{Review of Elementary Polymatroid Theory}

Recall that a function $f:2^{V} \rightarrow \mathbb{R}$ is submodular if $f(A) + f(B) \ge f(A \cup B) + f(A \cap B)$ for all $A,B \subseteq V$ and $g:2^{V} \rightarrow \mathbb{R}$ is supermodular if $g(A) + g(B) \le g(A \cup B) + g(A \cap B)$ for all $A,B \subseteq V$ . 

For the rest of this section we assume that  $f:2^{V} \rightarrow \mathbb{R}_+$ is a non-negative, non-decreasing (with respect to set inclusion) submodular function with $f(\emptyset)=0$, where $V=[n] \equiv \{1,\ldots,n\}$ for some positive integer $n$. (We set $[n]=\emptyset$ if $n=0$.) We assume that the values $f(S)$ are given by an oracle. 
Let $\mathcal{P}(f)$ be the polymatroid associated with $f$, given by
\[
\mathcal{P}(f) = \{x \in \mathbb{R}^n_+: x(S) \le f(S) \text{ for all } S \subseteq V\},
\]
where $x(S)\equiv \sum_{j \in S} x_j$. We first review the problem of maximizing a linear function $w^T x$ over $x \in \mathcal{P}(f)$, where $w \in \mathbb{R}^n_+$ is a constant. Let $\sigma:V \rightarrow V$ be a permutation (or bijection) of $V$ such that $w_{\sigma(1)} \ge \cdots \ge w_{\sigma(n)}$. The classic solution to the problem, given in \cite{edmonds1970submodular} is the point $x^\sigma$ given by
\begin{align}
x^\sigma_{\sigma(j)} = f(\{\sigma(1),\ldots,\sigma(j)\})-f(\{\sigma(1),\ldots,\sigma(j-1)\}),~ j=1,\ldots,n. \label{eq:vxs}
\end{align}
Notice that for any $w$ and $\sigma$, we have $x^\sigma(V)=f(V)$, so an equivalent problem is to maximize $w^T x$ over the {\em base polyhedron} $\mathcal{B}(f)$ of $f$, given by
\[
\mathcal{B}(f) = \{x \in \mathcal{P}(f): x(V)=f(V)\}.
\]
The vertices of $\mathcal{B}(f)$ are given by all points $x^\sigma$ defined by (\ref{eq:vxs}), as $\sigma$ ranges over the set $\Sigma \equiv \Sigma(V)$ of all possible permutations of $V$.

Now let $g$ be an arbitrary non-decreasing, {\em supermodular} function with $g(\emptyset)=0$. The contrapolymatroid $\mathcal{Q}(g)$ associated with $g$ is defined by
\[
\mathcal{Q}(g) \equiv \{x \in \mathbb{R}^n: x(S) \ge g(S) \text{ for all } S \subseteq V\}.
\]
The base of $\mathcal{Q}(g)$ is given by 
\begin{align}
\mathcal{B}(g) = \{x \in \mathbb{R}^n: x(S) \ge g(S) \text{ for all } S\subseteq V \text{ and } x(V)=g(V) \}. \label{eq:basepoly}
\end{align}
Vertices $x^\sigma$ of $\mathcal B(g)$ are given analogously to~(\ref{eq:vxs}).

Later, we will use the following fact, which is easy to verify.
\begin{lemma} \label{lem:prefix}
If $w_i > w_j$, then for any $x^\sigma$ that maximizes $w^T x$, there exists some $\tilde{\sigma}$ such that $x^{\tilde{\sigma}}=x^{\sigma}$ and $\tilde{\sigma}^{-1}(i) < \tilde{\sigma}^{-1}(j)$ (that is, $i$ precedes $j$ in $\tilde{\sigma}$).
\end{lemma}

We note that in giving running times, we assume that it takes only constant time to answer an oracle query.  

\subsection{The Main Problem}

The problem we consider in this paper is that of finding some $x \in \mathcal{B}(f)$ to maximize $\min_j w_j x_j$, where $f$ is an arbitrary non-decreasing submodular function with $f(\emptyset)=0$. This is a special case of the {\em lexicographically optimal base} problem, introduced by \cite{fujishige1980lexicographically}, where, subject to the minimum component being maximal, the second-smallest component is maximized, and so on.

The special case that we study is equivalent to a zero-sum game in which a pure strategy for Player~1 (the maximizer) is a permutation $\sigma$ of $V$ (or, equivalently, a vertex $x^\sigma$ of $\mathcal{B}(f)$) and a pure strategy for Player~2 (the minimizer) is a direction $j \in V$. For a given pair of pure strategies $\sigma$ and $j$, the payoff is given by
\[
P_{f,w}(\sigma,j) \equiv w_j x^\sigma_j.
\]
We will usually drop the $f$ and $w$ from the subscript of $P$. We denote this game by $\Gamma_{\max}(\mathcal{B}(f),w)$. We will also consider a variation of the game, which we denote by $\Gamma_{\max}(\mathcal{B}(f),w)$, which is identical except that Player~1 is the minimizer and Player~2 is the maximizer. Similarly, if $g$ is a non-decreasing supermodular function, we may consider the games $\Gamma_{\max}(\mathcal{B}(g),w)$ and $\Gamma_{\min}(\mathcal{B}(g),w)$, defined analogously.% where a pure strategy for Player~1 (the maximizer) is given by some vertex $x^\sigma$ of $\mathcal{B}(g)$, a pure strategy for Player~2 is some $j \in V$ and the payoff is $P_{\sigma,j}(x,y) \equiv w_j x^\sigma_j$. Interchanging the roles of minimizer and maximizer, we obtain our fourth and final variation of the game, which we denote by $\Gamma_{\max}(g,w)$.

 A mixed strategy for Player~1 in $\Gamma_{\max}(\mathcal{B}(f),w)$ corresponds to a point $x$ of $\mathcal{B}(f)$ and the expected payoff of such a strategy against a pure strategy $j$ of Player~2 is $w_j x_j$. 

A mixed strategy for Player~2 is a randomized choice of directions, where each $j \in V$ is chosen with some probability $\theta_j \ge 0$, where $\sum_{j=1}^n \theta_j = 1$. For such a mixed strategy, the payoff against a strategy $x$ of Player~1 is
\[
\sum_{j =1}^n \theta_j w_j x_j = x^T y,
\]
where $y = \sum_{j=1}^n \theta_j w_j e^j$, and $e^j$ is the $j$th coordinate vector. 

Equivalently, we may consider a mixed strategy for Player~2 as a point $y$ of the simplex
\begin{align*}
C = \big\{ \sum_{j=1}^n \theta_j w_j e^j : \sum_{j=1}^n \theta_j = 1 \text{ and } \theta_i \ge 0 \text{ for all } j=1,\ldots,n \big\},
\end{align*}
so that a pure strategy for Player~2 is a vertex $w_je^j$  of $C$. In a small abuse of our notation, we write $P(x,y)$ for the expected payoff $x^T y$ when Player~1 uses strategy $x$ and Player~2 uses strategy $y$. When one player uses a pure strategy and the other uses a mixed strategy, we extend the use of $P$ in the natural way.

Since each player has a finite number of pure strategies in each of its four versions, the game has optimal mixed strategies and a value $v$, by the minimax theorem for zero-sum games. For example, in the case of $\Gamma_{max}(\mathcal{B}(f),w)$,
\[
v=\max_{x \in \mathcal{B}(f)} \min_j P(x,j) = \min_{y \in C} \max_{\sigma \in \Sigma} P(\sigma,y).
\]

\section{Applications} \label{sec:apps}

In this section we show how our main problem and its variations can be used to model a number of search games as well as problems in sequential testing and queuing.

\subsection{Applications to Search Games} \label{sec:searchgames}

We begin by considering search games between a Searcher (Player~1) and a Hider (Player~2), where $V$ corresponds to a set of hiding locations. In each example, a Searcher pure strategy is a permutation $\sigma$ of $V$, where $\sigma(i)$ is the location that is in position~$i$ in the order of search and a Hider pure strategy is a location $i \in V$ at which a target is hidden.

\subsubsection{A weighted search game} \label{sec:weighted}
Consider a game where the time to search location $i$ is given by $t_i >0$ and each location $i$ has a weight $d_i$, corresponding to the rate of damage incurred at location $i$ while the target has not been found. The payoff is given by $P(\sigma, i)=d_i C_i^\sigma$, for a permutation $\sigma$ and $i \in V$, where
\[
C_i^\sigma = \sum_{\sigma^{-1}(j) \le \sigma^{-1}(i)} t_j.
\]
This payoff is the total time to find the Hider multiplied by the rate of damage. The Searcher is the minimizer and the Hider is the maximizer.  This game was considered by \cite{yolmeh-gursoy}, who solved the special case when the search times $t_i$ are all equal to~1, using a polyhedral approach. (\cite{yolmeh-gursoy} also applied a column and row generation approach to the game in a more general network setting, with multiple searchers and targets.)

\cite{condon2009algorithms} studied the special case of this game for $d_i=1/c_i$, which they called the {\em game theoretic multiplicative regret} game.  This case was also studied by \cite{angelopoulos2019expanding}. Implicit in the results of \cite{condon2009algorithms} is an optimal Player~1 (Searcher) strategy and the value of the game for the general weighted search game with arbitrary $d_i$.

Here we show how the game is a special case of $\Gamma_{\min}(\mathcal{B}(g),w)$. %give an alternative solution using the framework presented in Sections~\ref{sec:main} and~\ref{sec:min}. 
The searching of locations is analogous to the processing of jobs in single machine scheduling, and in the language of scheduling theory, we can interpret the time $t_i$ as the {\em processing time} of job $i$ and the time $C_i^\sigma$ as the {\em completion time} of job $i$ under the schedule $\sigma$. 
We associate a Searcher pure strategy $\sigma$ with a point $x^\sigma$ given by $x^\sigma_i = t_i C_i^\sigma,~ i \in V$.  
It is well known from scheduling theory (see \cite{queyranne1994polyhedral}) that the set of vectors $x^\sigma$ are the vertices of $\mathcal{B}(g)$, where $g$ is the supermodular function given by
\[
g(S) = \frac{1}{2}(t(S)^2 + t^2(S)),
\]
and $t^2(S) = \sum_{i \in S} t_i^2$. The polyhedron $\mathcal{B}(g)$ is known as the {\em scheduling polyhedron} and corresponds to the set of Searcher mixed strategies in the search game. Let $w_i = d_i/t_i$. Then for a Hider pure strategy $i$, the expected payoff against a Searcher mixed strategy given by~$x$ is $x_i w_i$. Hence, this is the game $\Gamma_{\min}(\mathcal{B}(g),w)$.

\subsubsection{A weighted search game with variable speeds} \label{sec:variable}

We can extend the model of the previous subsection by adopting the variable speed network model, as considering by \cite{alpern2014searching}. Suppose that we think of the set of locations $V$ as endpoints of $n$ arcs, whose other endpoint is a common point $O$. The Searcher successively travels from $O$ to the end of each arc and back again, where the time to travel from $O$ to the end of arc~$i$ is $a_i >0 $ and the time to travel back again is $b_i >0$. Let $t_i = a_i+b_i$ be the tour time of arc $i$. Similarly to the previous subsection, the vector $\tilde{C}^\sigma$ is defined by
\[
\tilde{C}_i^\sigma = a_i +  \sum_{\sigma^{-1}(j) < \sigma^{-1}(i)} t_j = C_i^\sigma - b_i,
\]
and corresponds to the times the Searcher reaches each location under $\sigma$. 

We consider a weighted search game with a minimizing Searcher and a maximizing Hider, whose payoff for a permutation $\sigma$ and $i\in V$ is given by $d_i \tilde{C}_i^\sigma$. If $b_i=0$ for all $i$, then $\tilde{C}^\sigma = C^\sigma$ and this is equivalent to the model of the previous subsection. 

The special case when the rates of damage  $d_i$ are all equal to~1 was solved by \cite{alpern2014searching} in the more general setting of tree networks, but the optimal Searcher strategy given had exponential support size even in the case of no network structure. The case of arbitrary $d_i$ has not been considered before. 

Let $\tilde{x}^\sigma_i = t_i \tilde{C}_i^\sigma$ and let $w_i = d_i/t_i$, so that the payoff for a Searcher strategy $\sigma$ and a Hider strategy~$i$ is $w_i \tilde{x}_i$. Note that we can write $\tilde{x}^\sigma = x^\sigma - c$, where $x^\sigma$ is defined as in the previous subsection and $c$ is given by $c_i=b_i t_i$. Therefore, the convex hull of the vectors $\tilde{x}^\sigma$ is equal to $\mathcal{B}(g)-c = \mathcal{B}(\tilde{g})$, where $\tilde{g}$ is the non-decreasing supermodular function given by
\[
\tilde{g}(S) = g(S)-c(S) = \frac{1}{2}(t(S)^2 +\sum_{j \in S} (a_j -b_j)t_j).
\]
Therefore, this is the game $\Gamma_{\min}(\mathcal{B}(\tilde{g}),w)$.

\subsubsection{A search and rescue game} \label{sec:s&r}

We now introduce a new search game in which we independently associate to every $i \in V$ a probability $p_i$ that the Searcher does not get captured when searching location $i$ and a probability $q_i$ that a target located at $i$ is found if location $i$ is searched. The payoff of the game is the probability the Searcher finds the target without getting captured herself. This is a generalization of the game introduced by \cite{lidbetter2020search} in which $q_i=1$ for all $i$.

More precisely, for a given permutation $\sigma$ and a given $i \in V$, the payoff is $q_i \pi_i^\sigma$, where
\[
\pi_i^\sigma = \prod_{\sigma^{-1}(j) \le \sigma^{-1}(i)} p_j. 
\]
The Searcher is the maximizer and the Hider is the minimizer. Let $x^\sigma_i = \frac{1-p_i}{p_i} \pi_i^\sigma$. It was shown by  \cite{kodialam2001throughput} and independently by \cite{agnetis2009sequencing} that the set of vectors $x^\sigma$ are the vertices of $\mathcal{B}(f)$ where $f$ is the non-decreasing submodular function given by
\begin{align}
f(S) = 1 - \prod_{i \in S} p_i. \label{eq:f}
\end{align}
Setting $w_i$ to be equal to $q_ip_i/(1-p_i)$, we see that this is the game $\Gamma_{\max}(\mathcal{B}(f),w)$.

\subsection{Relation to Sequential Testing} \label{sec:testing}

In this section we show that a sequential testing problem studied in \cite{condon2009algorithms} and \cite{kodialam2001throughput} is equivalent to the ``minimization'' version of the game considered in Subsection~\ref{sec:s&r}.  

Suppose some items, or {\em tuples} must be routed in some order through a set $V$ of operators, each of which tests whether the tuple satisfies some predicate (or filter) of a conjunction. To spread the load on the operators, different tuples may be routed in different orders. There is a known probability $p_i$ that a tuple will pass the test of operator $i$, and the tuple is routed through the operators until it fails one of the tests (and is eliminated) or it passes all of them. The problem here is to maximize the rate of flow of tuples routed through the operators, subject to the constraint that operator $i$ has a maximum flow rate of $r_i$. More precisely, the problem is given by the following linear program, where we denote the set of permutations of $V$ by $\Sigma(V)$.
\begin{align*}
\text{max } \sum_{\sigma \in \Sigma(V)} \lambda_\sigma \quad  \text{ s.t. } \quad \sum_{\sigma \in \Sigma(V)} \lambda_\sigma \prod_{\sigma^{-1}(j) < \sigma^{-1}(i)} p_j &\le r_i \text{ for all } i \in V, \\
\lambda_\sigma &\ge 0 \text{ for all } \sigma \in \Sigma(V).
\end{align*}
The variables $\lambda_\sigma$ here can be interpreted as the rate that tuples are routed through the operators in the order given by the permutation $\sigma$. We adopt the terminology of \cite{condon2009algorithms} and call this the {\em max-throughput} problem. The problem was solved in both \cite{condon2009algorithms} and \cite{kodialam2001throughput}, the latter paper exploiting the polymatroid structure of a space associated with the problem and the former giving a more efficient combinatorial algorithm with no reference to polymatroids.

Let $q_i = 1/(p_i r_i)$ and recall the notation $\pi_i^\sigma = \prod_{\sigma^{-1}(j) \le \sigma^{-1}(i)} p_j$ of the previous section. Let $v = 1/(\sum_{\sigma \in \Sigma(V)} \lambda_{\sigma})$ and let $\theta_\sigma = v \lambda_\sigma$. Then the max-throughput problem is equivalent to the following~LP
\begin{align*}
\text{min } v \text{ s.t. } \sum_{\sigma \in \Sigma(V)} \theta_\sigma q_i \pi_i^\sigma &\le v \text{ for all } i \in V, \\
\sum_{\sigma \in \Sigma(V)} \theta_\sigma &= 1, \\
\theta_\sigma &\ge 0 \text{ for all } \sigma \in \Sigma(V).
\end{align*}
This is the problem of finding an optimal strategy for Player~1 in the game $\Gamma_{\min}(\mathcal{B}(f),w)$, where $f$ is given by~(\ref{eq:f}) and $w_i=q_ip_i/(1-p_i)=r_i/(1-p_i)$. 

The derivation of the equivalence of these two problems closely follows the derivation in \cite{condon2009algorithms} of the equivalence of the game theoretic multiplicative regret problem and an artificial problem they called the {\em cumulative cost limit} problem.

\subsection{Applications to Queueing Theory} \label{sec:queueing}

As mentioned in the Introduction, the performance space of several multiclass queueing systems have been shown in \cite{federgruen1988characterization} and \cite{shanthikumar1992multiclass} to have a polymatroid structure. Possible performance measures of interest include the expected delay of the first $m$ jobs, the expected number of type $i$ jobs in the system at time $t$ or the expected number of job completions by time $t$. Depending on the context, the objective may be to maximize or minimize the performance measure and many such problems can be regarded as a special case of maximizing or minimizing a linear function over the base of a polymatroid.

For every maximization or minimization problem of this type we can consider a max-min or min-max variant. If we have an oracle for the submodular or supermodular function that defines the polymatroid or contrapolymatroid associated with a problem (in particular, if the function can be expressed in closed form), then the solution of the max-min or min-max problem follows from the results of this paper. We discuss one such problem here as an example rather than giving an exhaustive list of problems. 

\cite{coffman1980characterization} consider a queueing system with a single server with $n$ classes $V$ of jobs whose arrival times follow a Poisson process and whose service times are exponentially distributed (that is, a $M/M/1$ system). Jobs in class $i$ arrive at rate $\lambda_i$ and are serviced at rate $\mu_i$. The {\em traffic intensity} of jobs of class $i$ is $\rho_i=\lambda_i/\mu_i$. It is assumed that $\rho(V) \equiv \sum_{i=1}^n \rho_i <1$, which ensures the existence of a stationary distribution for the number of jobs in the system. The expected time that jobs of class $i$ spend in the system in the steady state is denoted $W_i$, and depends on the scheduling strategy chosen. 

Let $x \in \mathbb{R}^n$ be defined by $x_i=\rho_i W_i$. It is shown in \cite{coffman1980characterization} that the space of feasible vectors $x$ is the base $\mathcal{B}(g)$ of the contrapolymatroid given by the supermodular function
\[
g(S) = \frac{\sum_{i \in S} \rho_i/\mu_i}{1- \rho(S)}.
\]
Each vertex $x^\sigma$ of $\mathcal{B}(g)$ corresponds to a priority rule that assigns jobs to the server based on some fixed priority ordering of the job classes (given by the permutation $\sigma$). A non-vertex point $x = \sum_{\sigma \in \Sigma(V)} \theta_\sigma x^\sigma  \in \mathcal{B}(g)$ can be interpreted as a randomized priority rule where in each busy period the priority rule $\sigma$ is chosen with probability $\theta_\sigma$.

A well known consequence is that if the objective is to minimize some weighted sum $\sum_{i=1}^n c_i W_i$ of expected number of jobs in the system (where $c_i$ may correspond to the holding cost per unit time of jobs of class $i$), we can simply use the greedy algorithms of \cite{edmonds1970submodular} to minimize $w^T x$ with $w_i= c_i/\rho_i$. The solution is a priority rule that corresponds to some vertex of $\mathcal{B}(g)$.

Now suppose we wish to minimize the (weighted) {\em maximum} expected holding cost of any class of jobs. That is, we wish to find a performance vector $x \in \mathcal{B}(g)$ that minimizes $\max_i w_i x_i$. This is a special case of $\Gamma_{max}(\mathcal{B}(g),w)$.

\section{Solution and Special Cases} \label{sec:solution}

In this section we solve our main problem and its variations, and consider some special cases.  

\subsection{Solution to Main Problem}

We first note that for a given mixed strategy $y$ of Player~2 in the game $\Gamma_{\max}(\mathcal{B}(f),w)$, the problem of finding a best response for Player~1 is that of choosing $x \in \mathcal{B}(f)$ to maximize $x^T y$. This is the classical problem solved in \cite{edmonds1970submodular} of maximizing a linear function over $\mathcal{B}(f)$. 
With this observation, it follows that an optimal strategy for Player~1 can be computed in polynomial time (in $n$) using the ellipsoid algorithm
(see e.g.,~\cite{hellerstein2019solving}).
% LH I thought of including the following sentence but decided against it.  I'm leaving it here in case we want it later.
%The  
%problems of finding optimal strategies for Players~1 and 2 can be expressed as dual LPs, 
%and a separation oracle for the  Player~2 LP can be simulated using an algorithm solving the best-response problem for Player~1.  
\cite{fujishige1980lexicographically} showed that his (unique) solution to the lexicographically optimal base problem (and therefore an optimal Player 1 strategy) could be found in {\em strongly} polynomial time.  We give a new proof that this solution is an optimal Player 1 strategy. Our proof of optimality follows almost immediately from a duality approach.
%LH See above for mention of the ellipsoid algorithm.  I didn't say anything about ML approaches, which generally are not polynomial time.

For a subset $S \subseteq V$, $S \neq \emptyset$, denote $\sum_{i \in S} 1/w_i$ by $w^{-1}(S)$. Consider the Player~2 mixed strategy 
\[
y^S =  \sum_{i \in S} \left(\frac{w_i^{-1}}{w^{-1}(S)} \right) w_i e^i = \frac{1}{w^{-1}(S)} \sum_{i \in S} e^i. 
\]
For a Player~1 strategy $x \in \mathcal{B}(f)$, the expected payoff against $y^S$ is
\[
P(x, y^S) = \sum_{i \in S} x_i  \frac{1}{w^{-1}(S)} = \frac{x(S)}{w^{-1}(S)} \le \frac{f(S)}{w^{-1}(S)},
\]
by definition of $\mathcal{B}(f)$. We summarize this in the following lemma.
\begin{lemma} \label{lem:lb}
If Player~2 uses the strategy $y^S$ for some $S \neq \emptyset$, the expected payoff is at most $f(S)/w^{-1}(S)$.
\end{lemma}
We will show in Theorem~\ref{thm:value} that the strategy $y^S$ is optimal for Player~2, where $S$ is chosen to minimize $f(S)/w^{-1}(S)$. A minimizing set $S$ can be found in strongly polynomial time, using a parametric search (see \cite{iwata1997fast} [Section 6] for a parametric search algorithm for minimizing the ratio of a submodular function to a non-negative supermodular function). This relies on an algorithm for minimizing a submodular function. The fastest known strongly polynomial algorithm for submodular function minimization is that of \cite{orlin2009faster}, whose runtime is $O(n^6)$, so that the minimization of $f(S)/w^{-1}(S)$ takes time $O(n^7)$.

Before stating and proving the theorem, we define a strategy which will be optimal for Player~1. To do this, we recursively define a partition of $V$ into subsets $S_1,\ldots, S_r $. 

\begin{definition}[$f$-$w$ decomposition] \label{def:P1strat} Set $S_0 = \emptyset$ and suppose $S_0,\ldots,S_j$ have already been defined for some $j \ge 0$. Then if $S^j \equiv S_1 \cup \cdots \cup S_j$ is equal to $V$, set $r=j$. If not, we define $S_{j+1}$ to be any set $S \subseteq V \setminus S^j$ that minimizes $h_{S^j}(S)$, where
\[
h_{T}(S) \equiv \frac{f(T \cup S)-f(T)}{w^{-1}(S)}.
\]
We call $\mathcal{S}\equiv (S_1,\ldots,S_r)$ an {\em $f$-$w$ decomposition} of $V$. 
\end{definition}

Note that the function $h_{T}$ is the ratio of a submodular function and a modular function, therefore, as remarked earlier, it can be minimized in strongly polynomial time. Since $h_T$ is defined in terms of $f$ and $w$, a more informative notation is $h^{f,w}_T$, but we omit the superscripts in general when they are clear from the context.

We now define the Player~1 strategy $x^{\mathcal{S}}$ by
\[
x^{\mathcal{S}}_i = w_i^{-1} h_{S^{j-1}}(S_j) \text{ for all } i \in S_j, j=1,\ldots,r.
\] 
%Note that $x^{S^*}$ may not be uniquely defined. 
To show that $x^\mathcal{S}$ it is indeed a strategy, we need to prove that it lies in $\mathcal{B}(f)$.  Let $T \subseteq V$ be arbitrary and let $T_j = T \cap S_j$ for $j=0,1,\ldots,r$. Also set $T^{j} = \cup_{i \le j}T_i$. Then
\[
x^{\mathcal{S}}(T) = \sum_{j=1}^r \sum_{i \in T_j} w_i^{-1} h_{S^{j-1}}(S_j) = \sum_{j=1}^r w^{-1}(T_j) h_{S^{j-1}}(S_j) \le \sum_{j=1}^r f(S^{j-1} \cup T_j) - f(S^{j-1}),
\]
by definition of $S_j$. Since $f$ is submodular, $f(S^{j-1} \cup T_j) + f(T^{j-1}) \le  f(S^{j-1}) + f(T^j)$, so
\[
x^{\mathcal{S}}(T) \le \sum_{j=1}^r f(T^{j}) - f(T^{j-1}) = f(T).
\]
Hence, $x^{\mathcal{S}} \in \mathcal{P}(f)$. It is also easy to see that $x^{\mathcal{S}} (V)=f(V)$, so that $x^{\mathcal{S}}  \in \mathcal{B}(f)$.

It is elementary to show that the strategy $x^{\mathcal{S}}$ is actually the same for {\em any} $f$-$w$ decomposition, and is equivalent to Fujishige's solution to the lexicographically optimal base problem.

\begin{theorem} \label{thm:value}
%Let $f$ be a non-decreasing submodular function with $f(\emptyset)=0$ and 
Suppose $S^*$ is a non-empty set that minimizes $f(S)/w^{-1}(S)$. Then the value of the game $\Gamma_{\max}(\mathcal{B}(f),w)$ is equal to $f(S^*)/w^{-1}(S^*)$. An optimal strategy for Player~2 is $y^{S^*}$. An optimal strategy for Player~1 is~$x^{\mathcal{S}}$, where $\mathcal{S}=(S_1,\ldots,S_r)$ is any $f$-$w$ decomposition.
\end{theorem}
\textit{Proof.} 
By Lemma~\ref{lem:lb}, the value of the game is at most $f(S^*)/w^{-1}(S^*)$. To complete the proof, we will show that $x^\mathcal{S}$ ensures a payoff at least $f(S^*)/w^{-1}(S^*) = h_\emptyset(S_1)$ against any Player~2 strategy. Note that for a pure strategy $i$ of Player~2 with $i \in S_j$, the expected payoff against $x^\mathcal{S}$ is 
\[
P(x^\mathcal{S},i)=w_i x^{\mathcal{S}}_i  =  h_{S^{j-1}}(S_j).
\]
So it is sufficient to show that $h_{S^{j-1}}(S_j)$ is non-decreasing in $j$. By definition of $S_{j}$, we have
\begin{align}
\frac{f(S^j)-f(S^{j-1})}{w^{-1}(S_j)} \le \frac{f(S^{j+1})-f(S^{j-1})}{w^{-1}(S_j\cup S_{j+1})}, \label{eq:non-dec}
\end{align}
for $j=1,\ldots,r-1$. Writing $w^{-1}(S_j\cup S_{j+1}) = w^{-1}(S_j)+w^{-1}(S_{j+1})$ and rearranging yields
\[
w^{-1}(S_j)(f(S^{j+1})-f(S^j)) \ge w^{-1}(S_{j+1})(f(S^j)-f(S^{j-1})).
\]
This is equivalent to $h_{S^j}(S_{j+1}) \ge h_{S^{j-1}}(S_{j})$, and the proof is complete. 
\hfill $\Box$

Any given mixed strategy $y$ of Player~2 can be expressed uniquely as a convex combination of his pure strategies (that is, vertices of $w_j e^j$ of $C$) simply by taking $\theta_i = y_j/w_j$. A given mixed strategy $x$ of Player~1 can be written as a convex combination of at most $n$ of her pure strategies $x^\sigma$, by Carath\'eodory's Theorem. In general, as discussed in \cite{hoeksma2014decomposition}, such a representation can be found in strongly polynomial time by combining the generic approach of \cite{grotschel2012geometric} with the algorithm of \cite{fonlupt2009strongly} for finding the intersection of a line with a polymatroid. The runtime of this algorithm is $O(n^9)$. %As we will discuss in Section~\ref{sec:monotone},
For particular problems it is possible to exploit the structure of $\mathcal{B}(f)$ in order to find a more efficient algorithm for representing a Player~1 mixed strategy as a convex combination of at most $n$ of her pure strategies.

In general, both players have multiple optimal strategies. For Player~2, we can characterize these strategies. %The following lemma will be useful. Recall that the vertices of $\mathcal{B}(f)$ are the points $x^\sigma$, given in~(\ref{eq:vxs}) and $\Sigma$ is the set of all permutations of $V$.

Let $\mathcal{F}=\mathcal{F}(f)$ be the family of sets $S \neq \emptyset$ that minimize $f(S)/w^{-1}(S)$, so that the value~$v$ of the game is equal to $f(S)/w^{-1}(S)$ for any $S \in \mathcal{F}$. We also set $f(\emptyset)/w^{-1}(\emptyset)$ to be equal to $v$, so that $\emptyset \in \mathcal{F}$. It is useful to note that $\mathcal{F}$ is a lattice.  Indeed, suppose $S,T \in \mathcal{F}$. In the following calculation, we use the observation that for any $a,b,c,d >0$, if $a/b ,c/d \ge v$ then $(a+c)/(b+d) \ge v$, where the second inequality is tight if the first is also tight. We have 
\[
v = \frac{f(S)+f(T)}{w^{-1}(S) + w^{-1}(T)} \ge \frac{f(S \cup T) + f(S\cap T)}{w^{-1}(S \cup T) + w^{-1}(S \cap T)} \ge v,
\]
where the equality and second inequality follow from our observation and the first inequality follows from the submodularity of $f$. Therefore, the two inequalities hold with equality, and ${S \cup T, S \cap T \in \mathcal{F}}$.

\begin{theorem} \label{thm:uniqueness}
A Player~2 strategy $y$ is optimal if and only if it is in the convex hull of ${\{y^S: S \in \mathcal{F}(f)\}}$.
\end{theorem}
\textit{Proof.}
By Theorem~\ref{thm:value}, each element of $\{y^S: S\in \mathcal{F} \}$ is optimal, so any convex combination of such points is also optimal.  

For the opposite direction, suppose that $y^*$ is an optimal Player~2 strategy. By relabeling, let us assume that $y^*_1\ge \cdots \ge y^*_n$. Then recalling that $y^S = (\sum_{i \in S} e^i)/w^{-1}(S)$ for $S \subseteq V$ and setting $y^*_{n+1}=0$, we can write $y^*$ as
\[
y^* = \sum_{i=1}^n y^*_i e^i = \sum_{i=1}^n e^i \sum_{j =i}^n (y^*_j - y^*_{j+1}) = \sum_{j=1}^n (y^*_j - y^*_{j+1} )\sum_{i =1}^j e^i = \sum_{j=1}^n \lambda_j y^{[j]},
\]
where $\lambda_j = (y^*_j - y^*_{j+1})w^{-1}([j])$. Note that
\[
\sum_{j=1}^n \lambda_j = \sum_{j=1}^n (y^*_j - y^*_{j+1}) \sum_{i = 1}^j w_i^{-1}  = \sum_{i=1}^n w_i^{-1} \sum_{j = i}^n  (y^*_j - y^*_{j+1}) = \sum_{j=1}^n y^*_j/w_j = 1,
\]
where the final equality follows from the fact that $y^* \in C$. So $y^*$ is a convex combination of the strategies $y^{[j]}$. We claim that if $\lambda_k >0$ for some $k$ then $[k] \in \mathcal{F}$, so that $y^*$ is in fact a convex combination of strategies $y^S$ with $S \in \mathcal{F}$. Indeed, suppose that $\lambda_k >0$, so that $y^*_k > y^*_{k+1}$. Since any pure strategy best response $x$ to $y^*$ maximizes $x^T y^*$, by Lemma~\ref{lem:prefix}, we can express $x$ as a point $x^\sigma$ such that the first $k$ terms of $\sigma$ are $[k]$ in some order. So by definition of $x^\sigma$,
\begin{align}
\sum_{i =1}^k x_i = f([k]). \label{eq:sumx}
\end{align}
Equation~(\ref{eq:sumx}) also holds for any mixed strategy $x$ which is a best response to $y^*$ (since $x$ must be a mixture of pure best responses to $y^*$). In particular, it holds for $x=x^\mathcal{S}$, where $\mathcal{S}$ is any $f$-$w$ decomposition of $V$ whose first element $S_1$ is the maximal element $\cup_{S \in \mathcal{F}} S$ of $\mathcal{F}$. 

We claim that $[k] \subseteq S_1$. Let $i \in [k]$ and suppose $i \in S_j$ for some $j>1$. Since $y^*_i \ge y^*_k > y^*_{k+1} \ge 0$ and any Player~2 pure strategy in the support of $y^*$ that is played with positive probability must be a best response to $x^\mathcal{S}$, it follows that strategy $i$ is a best response to $x^\mathcal{S}$. But by the maximality of $S_1$, inequality~(\ref{eq:non-dec}) with $j=1$ is strict, and rearranging gives $h_\emptyset(S_1) < h_{S^{1}}(S_2)$. Since $h_{\emptyset(S^{j-1}}(S_j))$ is non-decreasing, for any $i' \in S_1$,
\[
P(x^\mathcal{S},i') = h_{\emptyset}(S_1) <  h_{S^{j-1}}(S_j) = P(x^\mathcal{S},i), 
\]
so $i$ cannot be a best response to $x^\mathcal{S}$, a contradiction. Hence, $i \in S_1$ so $[k] \subseteq S_1$.

Now, by definition of $x^\mathcal{S}$,
\[
\sum_{i =1}^k x^\mathcal{S}_i =\sum_{i=1}^k \frac{ w_i^{-1} f(S_1)}{w^{-1}(S_1)} = w^{-1}([k]) v,
\]
where $v$ is the value of the game. Combining this with~(\ref{eq:sumx}) yields $f([k])/w^{-1}([k]) = v$, so $[k] \in \mathcal{F}$. This completes the proof. 
\hfill $\Box$

\subsection{Special Cases}

To find optimal strategies in the game $\Gamma_{\max}(\mathcal{B}(f),w)$, it is necessary to minimize the function $h_T(S)=(f(T \cup S)-f(T))/w^{-1}(S)$. As previously remarked, there is a strongly polynomial time algorithm for this problem with runtime $O(n^7)$. To calculate an optimal Player~1 strategy, this algorithm must be run at most $n$ times, so the overall runtime is $O(n^8)$. For some functions $f$, this minimization can be performed much faster, as we show in the remainder of this section.

\begin{definition} We say that the payoff $P=P_{f,w}$ is {\em $\zeta$-decreasing} if there exists $\zeta \in \mathbb{R}^n_+$ such that for any $\sigma \in \Sigma(V)$ and any $i,j \in V$ with $\sigma^{-1}(i) < \sigma^{-1}(j)$,
\begin{align}
\frac{P(\sigma,i)}{P(\sigma,j)} \ge \frac{\zeta_i}{\zeta_j}. \label{eq:indexable}
\end{align}
If $\frac{P(\sigma,i)}{P(\sigma,j)} \le \frac{\zeta_i}{\zeta_j}$ we say $P$ is {\em $\zeta$-increasing}. If $\zeta_i=1$ for all $i$, then we say $P$ is {\em decreasing} (or respectively {\em increasing}). 
\end{definition}
If the payoff is $\zeta$-decreasing (or increasing) we assume that the values $\zeta_i$ are given as part of the input of the problem.

\begin{lemma} \label{lem:indexable}
Suppose $P=P_{f,w}$ is $\zeta$-decreasing. Then $S^* \equiv \cup_{S \in \mathcal{F}} S$ is equal to $\{i \in V: \zeta_i \le r\}$ for some $r >0$.
\end{lemma}

%\begin{lemma} \label{lem:indexable}
%Suppose there exists $z \in \mathbb{R}^n_+$ such that for any $S \subseteq V$ and any $i,j \notin S$,
%\begin{align}
%\frac{f(S \cup \{i\})-f(S)}{z_i} \ge \frac{f(S \cup \{i,j\}) - f(S \cup \{i\})}{z_j}. \label{eq:indexable}
%\end{align}
%Then $S^* \equiv \cup_{S \in \mathcal{F}} S$ is equal to $\{i \in V: w_iz_i \le r\}$ for some $r >0$.

%In particular, this holds for $z$-indexable functions $f$.
%\end{lemma}
\textit{Proof.} It is sufficient to show that if $P$ is $\zeta$-decreasing and $i \in S^*$ and $j \notin S^*$, then $\zeta_i <  \zeta_j$. Let $x^{\mathcal{S}}$ be any optimal Player~1 strategy such that the first set in the partition $\mathcal{S}$ is $S^*$, and write $x^{\mathcal{S}} = \sum_{\sigma \in \Sigma} \theta_\sigma x^\sigma$ as a convex combination of pure strategies. Since $y^{S^*}_i >  0 = y^{S^*}_j$, for any best response $\sigma$ to $y^{S^*}$, we can write $x^\sigma = x^{\tilde{\sigma}}$, where $\tilde{\sigma}^{-1}(i) < \tilde{\sigma}^{-1}(j)$, by Lemma~\ref{lem:prefix}. Since every pure strategy in the support of $x^{\mathcal{S}}$ must be a best response to $y^{S^*}$, we can assume that if $\theta_\sigma >0$ then $\sigma^{-1}(i) < \sigma^{-1}(j)$. It follows from~(\ref{eq:indexable}) that if $i \in S^*$ and $j \notin S^*$, then
\begin{align}
\frac{P(x^{\mathcal{S}},i)}{\zeta_i} =\sum_{\sigma \in \Sigma}   \frac{\theta_{\sigma} P(x^\sigma,i)}{\zeta_i} \ge   \sum_{\sigma \in \Sigma} \frac{\theta_{\sigma} P(x^\sigma,j)}{\zeta_j} = \frac{P(x^{\mathcal{S}},j)}{\zeta_j}. \label{eq:1}
\end{align}
By Theorem~\ref{thm:uniqueness}, every element of $S^*$ (in particular, $i$) is in the support of some optimal Player~2 strategy and $j$ cannot be in the support of any Player~2 strategy. Therefore, $i$ must be a best response to $x^{\mathcal{S}}$ and $j$ cannot be a best response, so that 
\begin{align}
P(x^{\mathcal{S}},i)  < P(x^{\mathcal{S}},j). \label{eq:2}
\end{align}
Combining~(\ref{eq:1}) and~(\ref{eq:2}) yields $\zeta_i <  \zeta_j$.
\hfill $\Box$

It is worth pointing out that although the definition of $\zeta$-decreasing and the proof of Lemma~\ref{lem:indexable} are given in game theoretic terms, the lemma is not exactly a game theoretic result, and could be stated without reference to the game $\Gamma(\mathcal{B}(f),w)$. Indeed, it is easy to see that $P_{f,w}$ is $\zeta$-decreasing if and only if there exists $\zeta' \in \mathbb{R}^n_+$ such that
\begin{align}
\frac{f(S \cup \{i\})-f(S)}{f(T \cup \{j\}) - f(T)} \ge \frac{\zeta'_i}{\zeta'_j}, \label{eq:equiv}
\end{align}
for any $S \subset T$ with $i \notin S,j \notin T$. 

Lemma~\ref{lem:indexable} implies that for games $\Gamma_{\max}(\mathcal{B}(f),w)$ with a $\zeta$-decreasing payoff function, the set $S^* = \cup_{S \in \mathcal{F}}S$ can be found in time $O(n \log n)$, simply by relabeling the the elements of $V$ so that they are in non-decreasing order of the index $\zeta_i$, computing $f([k])/w^{-1}([k])$ for each $k \in [n]$ and choosing the largest $k$ that minimizes this function. (Note that these $n$ computations can done in time $O(n)$ by keeping a record of $w^{-1}([k])$ each time and adding $w^{-1}_{k+1}$ to obtain $w^{-1}([k+1])$.) Therefore the value of the game $f(S^*)/w^{-1}(S^*)$ and the optimal Player~2 strategy $y^{S^*}$ can be found in time $O(n \log n)$. 

In order to compute the optimal Player~1 strategy $x^\mathcal{S}$ it is necessary to calculate an $f$-$w$ decomposition $\mathcal{S}$, which involves at most $n$ minimizations of functions of the form $h_T(S)$. It is easy to check that if $P(f,w)$ is $\zeta$-decreasing, then so is the function $P(f_T,w)$, where
\[
f_T(S) = f(T \cup S)-f(S).
\]
It follows that an $f$-$w$ decomposition can be found in time $O(n^2)$. (However, expressing $x^{\mathcal{S}}$ as a convex combination of at most $n$ pure strategies takes additional computation in general.)

We conclude this section by showing that when the payoff is decreasing, the solution of the game is particularly simple.

\begin{lemma} \label{lem:perm-dec}
%Let $f$ be a non-decreasing submodular function with $f(\emptyset) = 0$. 
If $P=P_{f,w}$ is decreasing then $f(S)/w^{-1}(S)$ is non-increasing in $S$ and the value of the game is $f(V)/w^{-1}(V)$. The strategy $x^{\mathcal{S}}$ is optimal for Player~1, where $\mathcal{S}$ consists only of the set $V$, and $y^{V}$ is optimal for Player~2.
\end{lemma}
\textit{Proof.} 
Let $S \neq \emptyset$ be a proper subset of $V$, and without loss of generality, assume that ${S=\{1,\ldots,k\}}$ for some $k$. Let $j \notin S$ and let $\sigma$ be any permutation of $V$ that starts with $(1,2,\ldots,k,j)$. Since $P$ is decreasing, for any $i \in S$,
\[
w_j (f(S \cup \{j\})-f(S)) = P(\sigma,j) \le P(\sigma,i) = w_i(f([i])-f([i-1])).
\]
Then setting $\theta_i = w_i^{-1}/w^{-1}(S)$, we obtain
\begin{align*}
\frac{f(S)}{w^{-1}(S)} &= \sum_{i=1}^k \theta_i w_i(f([i])-f([i-1]))\\
& \ge \sum_{i=1}^k \theta_i w_j (f(S \cup \{j\})-f(S)) \\
&= w_j (f(S \cup \{j\})-f(S)).
\end{align*}
Rearranging yields
\[
f(S \cup \{j\}) w^{-1}(S) \le f(S) w^{-1}(S \cup \{j\}),
\]
or equivalently,
\[
\frac{f(S \cup \{j\})}{w^{-1}(S \cup \{j\})} \le \frac{f(S)}{w^{-1}(S)}.
\]
This proves that $f(S)/w^{-1}(S)$ is non-increasing in $S$, so the value of the game is \\$\min_{S \subseteq V} f(S)/w^{-1}(S) = f(V)/w^{-1}(V)$.

The optimality of the stated strategies is immediate from Theorem~\ref{thm:value}.
\hfill $\Box$

\section{Other Variations of the Game} \label{sec:min}

Let $g^\#$ be the {\em dual} of $g$, given by $g^\#(S)=g(V)-g(V \setminus S)$. It is easy to show that $g^\#$ is submodular and non-decreasing with $g^\#(\emptyset)=0$ and $\mathcal{B}(g) = \mathcal{B}(g^\#)$.  Moreover, $P_{g,w}$ is $\zeta$-increasing if and only if $P_{g^\#,w}$ is $\zeta$-decreasing. Therefore, the game $\Gamma_{\max}(\mathcal{B}(g),w)$ is equivalent to $\Gamma_{\max}(\mathcal{B}(g^\#),w)$, and the solution follows immediately from Theorems~\ref{thm:value} and~\ref{thm:uniqueness}. Versions of Lemmas~\ref{lem:indexable} and~\ref{lem:perm-dec} also hold.

%\begin{theorem} \label{thm:super}
%Let $g$ be a non-decreasing supermodular function with ${g(\emptyset)=0}$ and suppose ${S^* \subseteq V}$ minimizes $g^\#(S)/w^{-1}(S)$. Then the value of the game $\Gamma_{\max}(g,w)$ is equal to ${g^\#(S^*)/w^{-1}(S^*)}$. A Player~2 strategy $y$ is optimal if and only if it is in the convex hull of $\{y^S: S \in \mathcal{F}(g^\#)\}$. The strategy $x^{\mathcal{S}}$ is optimal for Player~1, where $\mathcal{S}$ is any $g^\#$-$w$ decomposition of $V$.
%\end{theorem}

%So far, we have considered two equivalent zero-sum games where Player~1 is the maximizer and Player~2 is the minimizer. We now consider an alternative game with Player~1 as minimizer and Player~2 as maximizer. Let $g$ be a non-decreasing, supermodular function with $g(\emptyset)=0$, and let $\Gamma_{\min}(\mathcal{B}(g),w)$ be the game which is the same as $\Gamma_{\max}(\mathcal{B}(g),w)$ except that Player~1 is the minimizer and Player~2 is the maximizer. Similarly, for a non-decreasing, submodular function $f$ with $f(\emptyset)=0$, we define $\Gamma_{\min}(\mathcal{B}(f),w)$ to be the same as $\Gamma_{\max}(f,w)$ except that Player~1 is the minimizer and Player~2 is the maximizer.

The minimization version $\Gamma_{\min}(\mathcal{B}(f),w)$ of the game does not seem to be equivalent to the maximization version, but the solution and analysis are almost identical. We briefly describe the solutions here and leave the proofs as an exercise. 

Analogously to an $f$-$w$ decomposition for submodular $f$, for supermodular $g$ we define a $g$-$w$ max-decomposition $\mathcal{S}=(S_1,\ldots,S_r)$ as follows. Set $S_0 = \emptyset$ and suppose $S_0,\ldots,S_j$ have already been defined for some $j \ge 0$. Then if $S^j \equiv S_1 \cup \cdots \cup S_j$ is equal to $V$, set $r=j$. If not, we define $S_{j+1}$ to be any set $S \subseteq V \setminus S^j$ that {\em maximizes} $h^{g,w}_{S^j}(S)$. This time, the function $h^{g,w}_T$ is the ratio of a supermodular function and a modular function and can be maximized by using the procedure of \cite{iwata1997fast} to minimize the inverse ratio. Then the Player~1 strategy $x^{\mathcal{S}}$ is defined in precisely the same way as in the original version of the game. 

\begin{theorem} \label{thm:min} Let $f$ be a non-decreasing submodular function with $f(\emptyset)=0$ and let $g$ be a non-decreasing supermodular function with $g(\emptyset)=0$. Then the solutions to the games $\Gamma_{\max}(\mathcal{B}(f),w)$, $\Gamma_{\max}(\mathcal{B}(g),w)$, $\Gamma_{\min}(\mathcal{B}(g),w)$ and $\Gamma_{\min}(\mathcal{B}(g),w)$ are given in Table~\ref{tab:solutions}. The value and an optimal Player~1 strategy are indicated in the second and third columns of the table. In each case, the set of optimal Player~2 strategies is the convex hull of the set of $y^{S^*}$ where $S^*$ ranges over all possible values as given in the second column of the table. The fourth column gives a condition on the payoff for the set $S^*$ to have the form given in the fifth column. The sixth column gives a condition for $S^*$ to be equal to $V$.

\begin{table}[htb!]
\centering
\caption{Solutions to four versions of the game with submodular $f$ and supermodular $g$} 
\label{tab:solutions}
\begin{tabular}{|c|c|c|c|c|c|}
\hline
& & \textbf{$\mathcal{S}$ for optimal} & \textbf{Condition } & \textbf{$S^*$, if condi- } & \textbf{Condition} \\ 
\textbf{Game} & \textbf{Value} &  \textbf{Player~1}  & \textbf{on payoff} &   \textbf{tion on } & \textbf{on payoff}   \\ 
&& \textbf{strategy  $x^\mathcal{S}$} && \textbf{payoff holds}&  \textbf{for} $S^* = V$ \\ \hline
$\Gamma_{\max}(\mathcal{B}(f),w)$ & $\frac{f(S^*)}{w^{-1}(S^*)} = $ &$f$-$w$ &  $\zeta$-decreasing & $\{i: \zeta_i \le r\}$ & decreasing \\ 
&   $\min_{S \subseteq V} \frac{f(S)}{w^{-1}(S)}$ & decomposition &&& \\ \hline
$\Gamma_{\max}(\mathcal{B}(g),w)$ & $\frac{g^\#(S^*)}{w^{-1}(S^*)} = $ & $g^\#$-$w$  &  $\zeta$-increasing & $\{i: \zeta_i \le r\}$ &  increasing \\ 
&$\min_{S \subseteq V} \frac{g^\#(S)}{w^{-1}(S)}$ & decomposition& && \\ \hline
$\Gamma_{\min}(\mathcal{B}(g),w)$ &  $\frac{g(S^*)}{w^{-1}(S^*)} =$ & $g$-$w$ max- & $\zeta$-increasing & $\{i: \zeta_i \ge r\}$ &  increasing \\ 
& $ \max_{S \subseteq V} \frac{g(S)}{w^{-1}(S)}$ & decomposition&&& \\ \hline
$\Gamma_{\min}(\mathcal{B}(f),w)$ &  $\frac{f^\#(S^*)}{w^{-1}(S^*)} = $ & $f^\#$-$w$ max- & $\zeta$-decreasing & $\{i: \zeta_i \ge r\}$ & decreasing \\ 
&  $\max_{S \subseteq V} \frac{f^\#(S)}{w^{-1}(S)}$ &decomposition  &&& \\ \hline
\end{tabular}
\end{table}

\end{theorem}

\section{Implications for our Applications} \label{sec:impl}

We now discuss the implication of our results for the applications described in Section~\ref{sec:apps}.

\subsection{Weighted Search Games}

The solution of the weighted search game described in Subsection~\ref{sec:weighted} follows from Theorem~\ref{thm:min}. The value of the game is
\[
\max_{S \subseteq V} \frac{g(S)}{w^{-1}(S)} = \max_{S \subseteq V} \frac{(t(S)^2 + t^2(S))/2}{\sum_{i \in S} t_i/d_i}.
\]
It is easy to see that the payoff $P(\sigma, i)$ is $\zeta$-increasing where $\zeta=d$. Hence, by Theorem~\ref{thm:min}, the value and optimal strategies can be found in time $O(n \log n)$. To express the optimal Searcher strategy as a mixture of at most $n$ pure strategies, one can use the strongly polynomial time decomposition algorithm of \cite{hoeksma2014decomposition}.

We note that two different solutions of the special case when the rates of damage~$d_i$ are all equal to 1 were  given by \cite{lidbetter13} and \cite{alpern2013mining}, though in each solution the size of the support of the optimal Searcher strategy was exponential in $n$. \cite{condon2009algorithms} also considered this special case, calling it the {\em game theoretic total cost} problem. They found an optimal Searcher strategy of support size~$n$. Theorem~\ref{thm:min} implies an alternative polynomial time algorithm for finding an optimal Searcher strategy with support size~$n$. Furthermore, the payoff is increasing in this case, so Theorem~\ref{thm:min} implies that the optimal Hider strategy given by \cite{lidbetter13} and \cite{alpern2013mining} is unique.

The solution of the more general weighted search game with variable speeds of Subsection~\ref{sec:variable} also follows from Theorem~\ref{thm:min}. The value of the game is
\[
\max_{S \subseteq V} \frac{g(S)}{w^{-1}(S)} = \max_{S \subseteq V} \frac{(t(S)^2 + \sum_{j \in S} (a_j -b_j)t_j)/2}{\sum_{i \in S} t_i/d_i}.
\]
Again, the payoff function here is $\zeta$-increasing for $\zeta=d$, so the value and optimal strategies can be found in time $O(n \log n)$. Also, since $\mathcal{B}(\tilde{g})$ is simply a translation of $\mathcal{B}(g)$ by $-c$, we can again use the decomposition theorem of \cite{hoeksma2014decomposition} for $\mathcal{B}(g)$ to write an optimal mixed Searcher strategy $x \in \mathcal{B}(\tilde{g})$ as a convex combination of at most $n$ pure strategies.

For the special case considered in \cite{alpern2014searching} where  $d_i=1$ for all $i$, our solution here improves upon the optimal Searcher strategy of exponential support size. Also, since the payoff is increasing, Theorem~\ref{thm:min} implies that the optimal Hider strategy is unique.

\subsection{The Search and Rescue Game}

By Theorem~\ref{thm:value}, the value of the search and rescue game of Subsection~\ref{sec:s&r} is
\[
\min_{S \subseteq V} \frac{f(S)}{w^{-1}(S)} = \min_{S \subseteq V} \frac{1 - \prod_{i \in S} p_i}{\sum_{i \in S} (1-p_i)/(q_ip_i)}.
\]
The payoff is easily seen to be $\zeta$-decreasing where $\zeta_i=q_i$ (or indeed where $\zeta_i = q_i/p_i$). It follows from Theorem~\ref{thm:min} that the value and optimal strategies can be found in time $O(n \log n)$.

\cite{kodialam2001throughput} gave a strongly polynomial algorithm with runtime $O(n^3 \log n)$ for representing a point in $\mathcal{B}(f)$ as a convex combination of at most $n$ vertices, and we can use this to express the optimal Searcher strategy as a mixture of at most $n$ pure strategies.

In the special case considered by \cite{lidbetter2020search} where $q_i=1$ for all $i$, a solution was given but the size of the support of the optimal Searcher strategy was exponential in $n$. This approach gives an optimal strategy with support size $n$. Since the payoff is decreasing in this case, the optimal Hider strategy given in \cite{lidbetter2020search} is unique.

\subsection{Sequential Testing}

The solution to the sequential testing problem of Subsection~\ref{sec:testing} follows from Theorem~\ref{thm:min}.  The algorithm of \cite{kodialam2001throughput} is essentially a special case of the algorithm given in the proof of Theorem~\ref{thm:value}.

\subsection{Queuing Theory}

A solution to the queuing problem posed in Subsection~\ref{sec:queueing} is given by Theorem~\ref{thm:min} of this paper, and the value of the min-max expected holding cost is
\[
\max_{S \subseteq V} \frac{g(S)}{w^{-1}(S)} = \max_{S \subseteq V} \frac{\sum_{i \in S} \rho_i/\mu_i}{(1- \rho(S))\sum_{i \in S} \rho_i/c_i}.
\]

\section{Finding Optimal Strategies when the Payoff is Monotone} \label{sec:monotone}

%LH I did a lot of revising in this section to handle
% the implicit representation issue.  I removed any mention
% of an implicit representation, since I didn't think it was
% really relevant here, despite what I may have said before.
As mentioned in Section~\ref{sec:main}, expressing an optimal Player~1 strategy $x^S$ as a convex combination of pure strategies relies on an algorithm whose runtime is $O(n^9)$, in general. We have also seen that for particular polymatroids, this runtime can be reduced. In this section we show that if $f$ is submodular and $P_{f,w}$ is decreasing, then an optimal Player~1 strategy can be efficiently implemented.  More particularly, we show that
a random pure strategy can be drawn from the (exponentially-sized) support of this optimal strategy, with appropriate probability, in time $O(n)$.

\begin{theorem}
Suppose $f$ is submodular and $P_{f,w}$ is decreasing. Then there is an optimal Player~1 strategy $x$ for $\Gamma_{\max}(\mathcal{B}(f),w)$ such that a random pure strategy $x^{\sigma}$ for Player~1, drawn from the distribution on pure strategies defined by $x$, can be generated in $O(n)$ time.  An analogous result holds for $\Gamma_{\min}(\mathcal{B}(g),w)$ if $g$ is supermodular and $P_{g,w}$ is increasing.
\end{theorem}
\textit{Proof.} 
First, we introduce some notation. For $A \subseteq V$, let $f_A:2^{V \setminus A} \rightarrow \mathbb{R}_+$ be given by $f_A(S) = f(S \cup A)-f(A)$. We also write $f|_A$ for the function $f$ restricted to subsets of $A$ and $w|_A$ for the vector $w$  restricted to elements in $A$. It is easy to show that $f_A$ and $f|_A$ are submodular and the payoffs $P_{f|_A,w|_A}$ and $P_{f_A,w|_A}$ are decreasing.

We begin by constructing an optimal strategy for $P_{f,w}$.
We construct the strategy recursively.  If $n=1$, only one strategy is available, which is optimal. Suppose $n \ge 2$ and we have a construction for games such that the number of Player~2 strategies is $n-1$ and let $V'=V \setminus \{n\}$. Define
\begin{enumerate}
\item[(i)] $\Gamma_1 \equiv \Gamma_{\min}(\mathcal{B}(f|_{V'}),w|_{V'})$,
\item[(ii)] $\Gamma_2 \equiv \Gamma_{\min}(\mathcal{B}(f_{\{n\}}),w|_{V'})$,
\end{enumerate}
whose values are $V_1 \equiv f(V')/w^{-1}(V')$ and 
$V_2 \equiv (f(V)-f(\{n\}))/w^{-1}(V')$,
respectively, by Lemma~\ref{lem:perm-dec}. By induction, we have a construction for an optimal strategy for both of these games. Denote these optimal strategies $x^1$ and $x^2$, respectively. We now define two new strategies $\tilde{x}^1$ and $\tilde{x}^2$ for $\Gamma_{\min}(\mathcal{B}(f),w)$ as follows. The strategy $\tilde{x}^1$ is obtained by replacing each pure strategy $x^\sigma$ in $x^1$ with $x^{\sigma'}$, where $\sigma'$ is $\sigma$ followed by element $n$. The strategy $\tilde{x}^2$ is obtained by replacing each pure strategy $x^\sigma$ in $x^2$ with $x^{\sigma''}$, where $\sigma''$ is $\sigma$ preceded by element $n$.

Table~\ref{tab:payoffs} displays the payoff of the strategies $\tilde{x}^1$ and $\tilde{x}^2$ against the element $n$ and against any element of $V'$.

\begin{table}[htb!]
\centering
\caption{Expected payoffs $P(\tilde{x}^1,i)$ and $P(\tilde{x}^2,i)$ for $i=n$ and $i \in V'$.}
\label{tab:payoffs}
\begin{tabular}{c|c|c}
& $i=n$   & $i \in V'$ \\ \hline
$\tilde{x}^1$  & $w_n(f(V)-f(V'))$ &  $V_1$  \\ \hline
$\tilde{x}^2$  & $w_n f(\{n\})$  & $V_2$
\end{tabular}
\end{table}
The function $f(S)/w^{-1}(S)$ is non-increasing, by Lemma~\ref{lem:perm-dec}. Hence,
\[
w_n f(\{n\}) = \frac{f(\{n\})}{w^{-1}_n} \ge  \frac{f(V)}{w^{-1}(V)} = \frac{f(V)}{w^{-1}(V')+w^{-1}_n}. 
\]
Rearranging, we obtain
\[
w_n f(\{n\}) \ge  \frac{f(V)-f(\{n\})}{w^{-1}(V')} = V_2.
\]
Also, 
\[
\frac{f(V')}{w^{-1}(V')} \ge \frac{f(V)}{w^{-1}(V)} = \frac{f(V)}{1/w_n+w^{-1}(V')}.
\]
Rearranging gives
\[
w_n(f(V)-f(V'))  \le f(V')/w^{-1}(V') = V_1.
\]
It follows that by mixing appropriately between strategies $\tilde{x}^1$ and $\tilde{x}^2$, Player~1 can construct a strategy $x$ whose expected payoff against against any pure strategy (and therefore also any mixed strategy) of Player~2 is equal to some constant $c$. Therefore, by definition of the optimal Player~2 strategy, $y^V$,
\[
c=P(x, y^V) = \frac{f(V)}{w^{-1}(V)},
\]
so $c$ is the value $f(v)/w^{-1}(V)$ of the game and $x$ is optimal.

We note that, because each recursive call mixes between two strategies, the support of the final constructed strategy $x$ has size $2^n$.

We now describe how to generate a random
pure strategy $x^{\sigma}$ from the distribution
on pure strategies defined by $x$, without actually constructing $x$.  
The procedure is similar to the recursive construction above.  However, in each recursive call, we do not recursively
generate optimal strategies for both $\Gamma_1$ and $\Gamma_2$.
Instead, we first generate the payoffs in Table~\ref{tab:payoffs} and
calculate the mixing probabilities for $x^1$ and $x^2$, call them
$p'$ and $p''$ (=$1-p'$).
We then randomly choose between recursively generating a pure
strategy for $\Gamma_1$ or for $\Gamma_2$,
choosing the first with probability $p'$ and the second with probability
$p''$.  
Denote by $x^{\sigma}$ the pure
strategy that is generated.
If it was generated for $\Gamma_1$, we return $x^{\sigma'}$,
where $\sigma'$ is produced from $\sigma$ by appending element $n$.
If it was generated for $\Gamma_2$, we return $x^{\sigma''}$
where $\sigma''$ is produced from $\sigma$ by prepending element $n$.
It is clear that this procedure generates a random pure
strategy with the appropriate probability.

It remains to verify that this procedure can be implemented to run in time
$O(n)$.  Recall that we assume that each oracle query
can be answered in constant time.  The procedure makes $O(n)$ recursive calls.
The only non-trivial part of the analysis is the computation of
the mixing probailities $p'$ and $p''$ in a recursive call.  
These are computed from the four entries in Table~\ref{tab:payoffs}.
The entries in the first column of the table can be computed in
constant time.  The entries in the second column, $V_1$ and $V_2$,
are equal to $f(V')/w^{-1}(V')$ and $f(V)/w^{-1}(V')$ respectively.
Computing these values from scratch in each recursive call would take
linear time per recursive call.   
However, using the fact that $w^{-1}(V) = 1/w_n + w^{-1}(V')$, we can easily reduce the
computation in each recursive call to take constant time, 
by taking advantage of the computation done in the previous recursive call.
Thus the runtime is $O(n)$.

An analogous result for $\Gamma_{\max}(\mathcal{B}(g),w)$ can be proved similarly.
\hfill $\Box$

\section{Conclusion}

We have shown that a number of natural games that arise in different research areas can be understood and analyzed through a single unifying framework, allowing us to gain new insight into existing results and to prove new results.
%We have provided a unifying framework under which we can understand and analyze a number of natural games that arise in different research areas, gaining new insight into existing results and proving new results.
%We have defined a very general max-min problem, or equivalently a zero-sum game, and shown that its value and optimal strategies can be found in strongly polynomial time. We gave examples of several existing and new problems that fall into the framework of our game and its variations. 
There are many related problems in search theory and sequential testing that do not fall under this framework, including problems involving networks and multiple targets. A promising avenue for future research could be to explore polyhedral approaches to such problems.

\subsection*{Acknowledgements}
This material is based upon work supported by the National Science Foundation under Grant Numbers IIS-1909335 and IIS-1909446.

\bibliographystyle{apalike}
\bibliography{references}

\end{document}